\documentclass[12pt]{article}

\usepackage{authblk}
\usepackage{amsmath, amsthm, amssymb, epsfig}
\usepackage{mathtools}
\usepackage{mathrsfs}
\usepackage[dvipsnames]{xcolor}
\usepackage{color}
\usepackage{url}
\usepackage{algorithmicx}
\usepackage[ruled]{algorithm}
\usepackage{algpseudocode}
\usepackage{graphicx}
\usepackage{caption}
\usepackage{subcaption}
\usepackage{bbold}
\usepackage{fullpage} 

\usepackage{wrapfig}

\usepackage[scaled]{helvet} 
\usepackage{bm}

\newtheorem{theorem}{Theorem}
\newtheorem{lemma}{Lemma}
\newtheorem{corollary}[theorem]{Corollary}

\newtheorem{remark}{Remark}

\DeclareMathOperator*{\Null}{null}

\DeclareMathOperator*{\argmin}{arg\,min}

\newcommand{\vct}[1]{\bm{#1}}
\newcommand{\mtx}[1]{\bm{#1}}

\def \C {\mathbb{C}}

\def \E { \mathbb{E}}
\def \A {\mtx{A}}
\def \x {\vct{x}}
\def \y {\vct{y}}
\def \a {\vct{a}}
\def \amiss {\tilde{\A}}
\def \xopt {\vct{x}_\star}
\def \xcur {\vct{x}^{k+1}}
\def \xprev {\vct{x}_{k}}

\def \b {\vct{b}}
\def \br {\vct{r}}
\def \D {\mtx{D}}
\def \M {\mtx{M}}
\def \Amiss {\tilde{\mtx{A}}}

\def \W {\mathcal{W}}
 
\def \ai {\A_i}

\def \diag {{\rm diag }}
\def \proj {\mathcal{P}_{\mathcal{W}} }

\newcommand{\norm}[1]{\left\lVert#1\right\rVert}

\begin{document}

\title{Stochastic Gradient Descent for Linear Systems with \\Missing Data}

\author[1]{Anna Ma}
\author[2]{Deanna Needell}
\affil[1]{Institute of Mathematical Science, Claremont Graduate University}
\affil[2]{Department of Mathematics, University of California Los Angeles}

\maketitle
\begin{abstract}
Traditional methods for solving linear systems have quickly become impractical due to an increase in the size of available data. Utilizing massive amounts of data is further complicated when the data is incomplete or has missing entries. In this work, we address the obstacles presented when working with large data and incomplete data simultaneously. In particular, we propose to adapt the Stochastic Gradient Descent method to address missing data in linear systems. Our proposed algorithm, the Stochastic Gradient Descent for Missing Data method (mSGD), is introduced and theoretical convergence guarantees are provided. In addition, we include numerical experiments on simulated and real world data that demonstrate the usefulness of our method. 
\end{abstract}

\section{Introduction}
When handling large amounts of data, it may not be possible to load the entire matrix (data set) into memory, as typically required by matrix inversions or matrix factorization. This has led to the study and advancement of stochastic iterative methods with low memory footprints such as Stochastic Gradient Descent, Randomized Kaczmarz, and Randomized Gauss-Seidel~\cite{ma2015convergence, needell2014stochastic, strohmer2009randomized, leventhal2010randomized}. The need for algorithms that can process large amounts of information is further complicated by incomplete or missing data, which can arise due to, for example, attrition, errors in data recording, or cost of data acquisition. Standard methods for treating missing data, which include data imputation~\cite{efron1994missing, fichman2003multiple}, matrix completion~\cite{cai2010singular, keshavan2009matrix, recht2011simpler, keshavan2010matrix}, and maximum likelihood estimation~\cite{dempster1977maximum, little2014statistical} can be wasteful, create biases, or be impractical for extremely large amounts of data. This work simultaneously addresses both issues of large-scale and missing data.

Consider the system of linear equations $\A\x = \b$\footnote{The linear system is not assumed to be consistent; we will use the notation $\A \x=\b$ to denote a general linear system. }, where $\A \in\C^{m\times n}$ is a large, full-rank, overdetermined ($m > n$) matrix. Suppose that $\A$ is not known entirely, but instead only some of its entries are available. As a concrete example, suppose $\A$ is the rating matrix from the survey of $m$ users about $n$ service questions, and $\b$ contains the $m$ ``overall'' ratings from each user (which is fully known). Each user may not answer all of the individual service questions, but a company wishes to understand how each question affects the overall rating of the user. That is, given partial knowledge of $\A$, one wishes to uncover $\xopt = \argmin_{\x} \frac{1}{2m} \| \A \x - \b  \|^2$. 

Let $\Amiss = \D \circ \A$ where $\A$ denotes the full matrix, and $\circ$ be the element-wise product, $\D$ denotes a binary matrix (1 indicating the availability of an element and 0 indicating a missing entry). Formally, one wants to solve the following optimization program:
\begin{equation}
\begin{split}
&\text{Given } \Amiss, \, \b \text{ s.t. }   \A \x = \b  \,\,\, and\,\,\, \Amiss = \D \circ \A  \\
&\text{Find } \xopt = \argmin_{\x \in \mathcal{W}} \frac{1}{2m}\| \A \x - \b  \|^2,
\end{split}
\label{eq:problem}
\end{equation} 
where $\mathcal{W}$ is a convex domain containing the solution $\xopt$ (e.g. a ball with large enough radius). 

\textbf{Contributions. } This work presents a stochastic iterative projection method for solving large-scale linear systems with missing data. We provide theoretical bounds for the proposed method's performance and demonstrate its usefulness on simulated and real world data sets.

\subsection{Stochastic Gradient Descent}
Stochastic iterative methods such as Randomized Kaczmarz (RK) and Stochastic Gradient Descent (SGD) have gained interest in recent years due to their simplicity and ability to handle large-scale systems. Originally discussed in~\cite{robbins1951stochastic}, SGD has proved to be particularly popular in machine learning~\cite{bottou2010large, zhang2004solving, bottou2012stochastic}. SGD minimizes an objective function $F(\x)$ over a convex domain $\mathcal{W}$ using unbiased estimates for the gradient of the objective, i.e., using $f_i(\x)$ such that $\E [\nabla f_i(\x) ] = \nabla F(\x)$. At each iteration, a random unbiased estimate, $\nabla f_i(\x)$, is drawn and the minimizer of $F(\x) $ is estimated with:
\begin{equation}\label{eq:sgd}
\xcur = \proj \left(\xprev - \alpha_k \nabla f_i(\xprev) \right),
\end{equation}
where $\alpha_k$ is an appropriately chosen step size, or learning rate, at iteration $k$ and $\proj$ denotes the projection onto the convex set $\mathcal{W}$. To solve an overdetermined linear system $\A \x = \b$, one approach is to minimize the least-squares objective function $F(\x) = \frac{1}{2m} \| \A \x - \b \|^2 = \frac{1}{m} \sum_{i=1}^m  f_i(\x)$ where $f_i(\x) = \frac{1}{2} (\A_i \x - \b_i)^2$, $\A_i$ denotes the $i^{th}$ row of $\A$, and $\b_i$ denotes the $i^{th}$ entry of $\b$. In this setting, a random \emph{row} of the matrix $\A$ is selected and \eqref{eq:sgd} is computed with $\nabla f_i(\xprev) = \A_i^*(\A_i \xprev - \b_i) $ where $(\cdot)^*$ denotes the conjugate transpose. 

The performance of SGD on linear systems depends on the choice of $\alpha_k$ and the consistency of the system (i.e. whether a solution to the system exists). When the linear system is consistent, SGD achieves linear convergence with an appropriately chosen fixed step size~\cite{schmidt2013fast}. For example, RK, a special instance of SGD for linear systems, has been shown to converge linearly for consistent systems without decreasing step sizes~\cite{kaczmarz1937angenaherte, strohmer2009randomized, needell2014stochastic}. Unfortunately, this is not the case when the system is inconsistent. When the linear system is inconsistent, or $\A \x \approx \b$, one must use decreasing step sizes to obtain the optimum (see e.g.~\cite{schmidt2013fast, hanke1990acceleration, censor1983strong}). This phenomenon is explained by the norm of the unbiased estimates at the minimizer, $\| \nabla f_i (\xopt)\|^2$. For consistent systems, $ \| \nabla f_i(\xopt) \|^2 = \| \ai^*(\ai \xopt - \b_i) \|^2 = 0$ since $\A \xopt = \b$. Intuitively, as SGD progresses closer to the minimizer, the magnitude of the iterates get smaller and allow SGD to converge. When the system is inconsistent, $\A \xopt = \b + \br$ for some residual vector $\br$ and least squares minimizer $\xopt$.As the SGD approximates approach $\xopt$, the magnitude of the iterates do not converge to 0 since $ \|\nabla f_i(\xopt) \|^2 = \| \ai^*(\ai \xopt - \b_i) \|^2 = \br_i^2\| \ai \|^2$. Using diminishing step sizes dampens the magnitude of the iterates over time, allowing SGD to converge. When SGD with fixed step size is applied to inconsistent systems, the iterates oscillate within a fixed distance from the solution~\cite{needell2014stochastic}. The fixed distance, also referred to as the convergence horizon, is proportional to the step size but inversely proportional to the rate of convergence. Therefore, there is a trade-off between the rate of convergence (speed) and the radius of convergence (accuracy). 

The proposed method, which we refer to as mSGD, is an SGD-type iterate with a correction term that takes into account the fact that not all entries of $\A$ are available. We start with a discussion on the model under which mSGD operates and proceed to derive the iterate. After the introduction of the algorithm, the formal results are stated. 

\textbf{Outline. } Section~\ref{sec:msgd} introduces the proposed method, the Stochastic Gradient Descent for Missing Data method (mSGD) and the main theoretical results. The performance of mSGD on simulated and real world data are shown in Section~\ref{sec:exp}. Finally, we conclude in Section \ref{sec:conclu}.

\section{Stochastic Gradient Descent for Missing Data}
\label{sec:msgd}

We model whether an entry of $\A$ is missing with i.i.d. Bernoulli random variables that are equal to $1$ with probability $p$. Practically, there are many applications in which this type of assumption holds. For example, surveys where participants are given a random subset of questions to answer follow this assumption. In collaborative filtering, there are various models where such assumptions hold \cite{goldberg1992using, marlin2007collaborative}. As another example, consider an extremely large $m \times n$ matrix $\A$ where it is not possible to load entire rows of $\A$ nor columns of $\A$ due to memory constraints. Instead, one is restricted to only loading $\lfloor p n \rfloor $ (random) elements of $\A$ at a time. Under this probabilistic assumption on the missing entries, the least squares solution can be computed without making any additional assumptions on the structure of $\A$ such as sparsity or low-rankness. In the case of having a fixed matrix $\Amiss \in \mathbb{C}^{m \times n}$ with missing entries, the theoretical results hold only if each row of the matrix is utilized once. If $\Amiss$ is an extremely overdetermined matrix (i.e. $m \gg n$), then this is a reasonable assumption.

\textbf{Notation.} Let $\D$ be an $m \times n$ matrix where the entries of $\D$, denoted by $\delta_{i,j}$ for $i=1, 2, \ldots m$ and $j= 1, 2, \ldots n$, are drawn independent and identically distributed (i.i.d.) from a Bernoulli distribution with parameter $p$ so that $\delta_{i,j} = 1$ with probability $p$. The matrix $\D$ is referred to as a \textit{binary mask} throughout and its entries indicate the locations of non-missing entries of $\A$. Let $\D_i$ be the diagonal matrix whose diagonal is equal to the $i$th row of $\D$. Given an $n \times n$ matrix $\M$, we denote the a matrix containing only the diagonal of $\M$ as $\diag(\M)$. Let $\Amiss$ represent the matrix $\A$ with missing elements filled in with zeros so that $\Amiss = \D \circ \A$ and $\Amiss_i = \D_i \A_i^*$, where $\circ$ denotes the element-wise product. Additionally, let $\sigma_{\min}(\A)$ be the smallest singular value of $\A$ and $\| \cdot \|$ denote the $\ell_2$-norm. The expected value taken over the random selection of rows of $\tilde{\A}$ is denoted $\E_i [\cdot]$, the expected value taken over all ($2^{mn}$) possible binary masks $\D$ as $\E_\delta [\cdot]$, and the full expected value as $\E[\cdot]$. Lastly, let $\mathcal{W}$ be some convex domain containing $\xopt$ and $B := \max_{\x \in \mathcal{W}} \| \x \|^2$. 

\subsection{The method}
Suppose one naively applies SGD to the system $\Amiss \x = \b$. To that end, consider the objective $\widehat{F} = \frac{1}{2m} \| \Amiss \x - \b \|^2 = \frac{1}{m} \sum_{i=1}^m \widehat{f}_i(\x)$ where $\widehat{f}_i (\x) = \frac{1}{2} (\amiss_i \x - \b_i)^2 $. 
This objective function leads to the update:
\begin{equation*}
\xcur = \xprev - \alpha_k \left(\amiss_i^*(\amiss_i\xprev - \b_i)\right),
\end{equation*}
since $\nabla \widehat{f}_i (\x) = \amiss_i^*(\amiss_i\x - \b_i) $.
Unfortunately, one computes that, taking the expectation with respect to the binary mask and gradient direction, 
\begin{align*}
\E_i \E_\delta [\nabla \widehat{f}_i (\x)] &= \E_i \E_\delta [\amiss_i^*(\amiss_i\x - \b_i) ]  \\
&= \frac{1}{m} \left( p^2 \A^*\A\x + (p-p^2)\diag(\A^*\A)\x - p \sum_i \ai^* \b_i \right) \\
&\neq \nabla F(\x).
\end{align*}
As a result, the iterates are not moving in the gradient descent direction toward the desired solution in expectation. 

Now, since we have information on the distribution of missing entries, we can use this to design a better objective function. For example, we can approximate the proportion of the right hand side vector $\b$ which can be accounted for using the distribution for missing entries. In other words, since $\E_\delta [\Amiss_i \x] = p \b_i$, consider the objective $\tilde{F}(\x) = \| \Amiss \x - p \b \|^2_2 $. Applying SGD to this objective, one computes
\begin{align*}
\E_i \E_\delta [\nabla \tilde{f}_i (\x)] &= \E_i \E_\delta [\amiss_i^*(\amiss_i\x - p\b_i) ]  \\
&= \frac{1}{m} \left( p^2 \A^*\A\x + (p-p^2)\diag(\A^*\A)\x - p^2 \sum_i \ai^* \b_i \right) \\
&\neq \nabla F(\x),
\end{align*}
which is again not the direction that one wants on average.

Instead of using $\nabla \tilde{f}(\x)$ as the step direction, we use $\nabla \tilde{f}_i(\x)$ to estimate $\nabla F(\x)$. In other words, we want to represent $\nabla F(\x)$ in terms of $\E [\nabla \tilde{f}_i(\x)]$. By doing so, iterates $\xcur$ move in the gradient descent direction towards the least squares solution to the objective $F(\x) = \frac{1}{2m} \| \A \x - \b \|^2$. From the above computation, one can see
$$ 
\nabla F(\x) = \frac{1}{p^2} \E [\nabla \tilde{f}_i(\x)] - \frac{(1-p)}{p^2} \E[\diag(\amiss_i^* \amiss_i)] \x.
$$
The detailed computation is available in the Appendix (Lemma~\ref{lem:egx}). Therefore the appropriate update is 
$$
\xcur = \xprev - \alpha_k \left( \frac{1}{p^2} \left(\amiss_i^*(\amiss_i\xprev - p\b_i) \right) - \frac{1-p}{p^2}\diag(\amiss_i^*\amiss_i)\xprev \right).
$$
Note that in classical SGD literature, the expected value is taken over the row choice $i$ when being applied to linear systems. However, in this setting there are two sources of randomness: the randomness from row selection and the randomness incurred by modeling missing data. In this computation, the expected value is taken with respect to both sources of randomness. The method is outlined in Algorithm~\ref{alg:msgd}.

\begin{algorithm}
\begin{algorithmic}[1]
\Procedure{}{$\Amiss$, $\b$, $T$, $p$, $\{\alpha_k\}$} \Comment{If using a fixed step size $\alpha$, $\alpha_k = \alpha$ for all $k$.  }
\State Initialize $\x_0$
\For {$k=1,2,\ldots, T $} 
	\State Choose row $i$ of $\Amiss$ with probability $\frac{1}{m}$ 
   \State $g(\xprev) = \frac{1}{p^2} \left(\amiss_i^*(\amiss_i\xprev - p \b_i) \right) - \frac{1-p}{p^2}\diag(\amiss_i^*\amiss_i)\xprev$ \label{algline:g}
   \State $\xcur = \proj \left( \xprev - \alpha_k g(\xprev) \right) $ \Comment{$\proj$ is the projection onto the set $\mathcal{W}$.}
\EndFor
\State Output $\xcur$
\EndProcedure
\end{algorithmic}
\caption{Stochastic Gradient Descent for Missing Data (mSGD)}\label{alg:msgd}
\end{algorithm}

\subsection{Main Results}

Before the main results are presented, note the following properties of the objective function,
\begin{equation}
F(\x) = \frac{1}{2m} \| \A \x - \b \|^2,
\label{eq:obj}
\end{equation}
and the update function in Algorithm~\ref{alg:msgd} (Line~\ref{algline:g}),
\begin{equation}
g(\x) = \frac{1}{p^2} \left(\amiss_i^*(\amiss_i\x - p \b_i) \right) - \frac{(1-p)}{p^2}\diag(\amiss_i^*\amiss_i)\x,
\label{eq:update}
\end{equation}
as they play an important role in the convergence analysis of mSGD.
\begin{itemize}

\item The objective function~\eqref{eq:obj} is $\mu$-strongly convex. For all $\x, \y \in \mathcal{W}$, 
\begin{equation*}
(\x - \y)^*(\nabla F(\x) - \nabla F(\y) ) \geq \mu \|\x - \y \|^2,
\end{equation*}
where 
\begin{equation}
\mu = \frac{\sigma^2_{min}(\A)}{m}.
\label{eq:mu}
\end{equation}

\item The update function $g(\x)$ is Lipschitz continuous, has Lipschitz constant $L_{i,D}$ (for a fixed instance of $i$ and $D$), and supremum Lipschitz constant $L_g$. In other words, for all $\x, \y \in \mathcal{W}$,
\begin{align}
\| g(\x) - g(\y) \| &\leq L_{i,D} \| \x - \y \| \label{eq:lid}\\
L_g &= \sup_{i,D} L_{i,D}. \label{eq:lg}
\end{align}
The supremum is taken over all choices of rows and all possible binary masks (i.e. all $2^{mn}$ possible binary masks). 
 
\item There exists a constant $G$ that uniformly bounds the expected norm of $\| g(\x) \|^2$,
\begin{equation}
\E[ \| g(\x) \|^2 ] \leq G,
\label{eq:G} 
\end{equation}
for all $\x \in \mathcal{W}$ and rows $\Amiss_i$. The expected norm of $g(\xopt)$ plays an important role in the convergence horizon. For this reason, let $G_\star$ denote the upper bound of $\E [\| g(\xopt) \|^2 ]$:
\begin{equation}
\E [\| g(\xopt) \|^2 ] \leq G_\star.
\label{eq:Gstar}
\end{equation}
\end{itemize}
The computation of $L_{i,D}$ and $L_g$ are shown in Lemma~\ref{lem:lg}. Lemma~\ref{lem:G} shows the computation of $G$ and $G_\star$. The statements and proofs of both lemmas are provided in the Appendix so that we may proceed to the presentation of the main results.

Theorem~\ref{thm:update} shows that, in expectation, Algorithm~\ref{alg:msgd} converges to the least squares solution of the linear system $\A \x = \b$ with properly chosen step size. This theorem is an application of the previously proven result stated in Lemma~\ref{lem:sgdbound}. The fixed step size regime and the trade off between convergence rate and accuracy is explored in Theorem~\ref{thm:fixed}. In addition, we provide an optimal step size choice based on a desired error tolerance, $\epsilon$, and a bound on the number of iterations required to obtain said tolerance in Corollary~\ref{cor:fixed}. Lastly, we remark on the recovery of classical SGD when $p=1$ both algorithmically and with respect to the proven error bounds.

\begin{theorem} Consider~\eqref{eq:problem} with $\Amiss = \D \circ \A$ where entries of $\D$ are drawn i.i.d. from a Bernoulli distribution with probability parameter $p$. Let $\mu$ be as defined in~\eqref{eq:mu}.  Choosing $\alpha_k = \frac{1}{\mu k}$, Algorithm~\ref{alg:msgd} converges in expectation with error
\begin{equation*}
\E[ \| \xcur - \xopt \|^2] \leq \frac{17 G(1+\log(k)) }{\mu^2 k},
\end{equation*} 
\label{thm:update}
where $G = \frac{2B}{mp^2}\left( 1 + \frac{(1-p)(2 - p)}{p}  \right)  \sum_i \norm{ \A_i }^4  + \frac{2}{mp^2} \sum_i \norm{ \A_i }^2 | \b_i|^2$ is an upper bound on $\E[\|g(\x)\|^2]$ and $B = \max_{\x \in \mathcal{W}} \| \x \|^2$.
\end{theorem}

It is clear that the convergence behavior of Algorithm~\ref{alg:msgd} depends on $G$, the uniform upper bound on the expected norm of $g(\x)$ and on $\sigma^2_{min}(\A)$. As one would expect, the more data that is missing, the larger the upper bound on expected error. In particular, assuming all other variables are constant and $p \in (0, 1]$, as $p$ decreases, $G$ increases. Theorem~\ref{thm:update} is an application of the following previously proved lemma.

\begin{lemma}(\cite{shamir2013stochastic} Theorem 1) Let $F(\x)$ be a $\mu$-strongly convex objective function, $g(\x)$ be such that $\E [ g(\x) ] = \nabla F(\x)$, and $\E[ \| g(\x) \|^2 ]\leq G$ for all $\x \in \mathcal{W}$. Using step size $\alpha_k = \frac{1}{\mu k}$ and update $\xcur = \mathcal{P}_\W (\xprev - \alpha_k g(\xprev)) $, it holds that
$$\E [ F(\xcur) - F(\xopt) ] \leq \frac{17G(1 + \log(k))}{\mu k}. $$
\label{lem:sgdbound}
\end{lemma}

The next theorem details the convergence behavior of Algorithm~\ref{alg:msgd} when using a fixed step size. Theorem~\ref{thm:fixed} shows that Algorithm~\ref{alg:msgd} experiences a convergence horizon that depends on $L_g$ and $G_\star$. For $p \in (0, 1]$, as $p$ decreases, $G_\star$ and $L_g$ both increase. Intuitively this makes sense as a larger amount of missing data should increase the size of the convergence horizon. Additionally, the convergence rate $r = \left(1 - 2\alpha \mu \left(1- \alpha L_g\right) \right)$ also increases as $p$ decreases. In other words, more missing data causes a slower convergence rate.

\begin{theorem}  Consider~\eqref{eq:problem} with $\Amiss = \D \circ \A$ where entries of $\D$ are drawn i.i.d. from a Bernoulli distribution and are equal to 1 with probability $p$. Let $L_g$, $G_\star$, and $\mu$ be as defined in~\eqref{eq:lg},~\eqref{eq:Gstar}, and~\eqref{eq:mu} respectively. Additionally, let the fixed step size be $\alpha < \frac{1}{L_g}$. Algorithm~\ref{alg:msgd} converges with expected error
\begin{equation}
\E[ \| \xcur - \xopt \|^2 ]\leq r^k \| \x_0 - \xopt \|^2 + \frac{\alpha G_\star}{\mu \left(1- \alpha L_g\right)}.
\label{eq:main}
\end{equation}
where $r = \left(1 - 2\alpha \mu \left(1- \alpha L_g\right) \right)$, $L_g = \frac{1}{p^2} sup_{i} \| \ai \|^2$, and $\mu = \frac{\sigma^2_{min}(\A)}{m}$. If $\A \x = \b$ is consistent, $G_\star = \frac{2(1-p)(2 - p)}{mp^3}  \| \xopt \|^2 \sum_i \norm{ \A_i }^4$. If the linear system is inconsistent (i.e. $\A \x = \b + \br$ for some residual vector $\br$), then $G_\star = \frac{2}{mp^2} \sum_i \norm{ \A_i }^2  \br_i^2  + \frac{2(1-p)(2 - p)}{mp^3}  \| \xopt \|^2 \sum_i \norm{ \A_i }^4$.
\label{thm:fixed}
\end{theorem}

Corollary~\ref{cor:fixed} and the subsequent remark comment on the number of iterations required by Algorithm~\ref{alg:msgd} to obtain some desired error tolerance $\epsilon$ using a particular fixed step size $\alpha^*$. The corollary itself details this information in terms of the variables in Theorem~\ref{thm:fixed} while the remark translates and simplifies $\alpha^*$ and $k$ (the number of iterations) into terms relating to $\A$. Note that the number of iterations required to reach a specified tolerance is a function of the ratio between the log of the initial error and $\epsilon$. The number of iterations increase as $\epsilon$ decreases. Additionally, the remark shows that as $p$ decreases, or as less data becomes available, more iterations are required to obtain an expected error of $\epsilon$. The proof of Corollary~\ref{cor:fixed} can be found in (\cite{needell2014stochastic} Corollary 2.2) with different constants. 

\begin{corollary} Given an initial error $\epsilon_0$ and choosing the fixed step size 
$$ \alpha^* = \frac{\epsilon \mu }{2G_\star + 2\mu \epsilon L_g }, $$
after 
$$ k = 2 \log \left( \frac{2 \epsilon_0}{\epsilon} \right) \left(\frac{L_g}{\mu} + \frac{G_\star}{\mu^2 \epsilon}\right)$$
iterations of Algorithm~\ref{alg:msgd}, $\E[ \| \xcur - \xopt \|^2 ] \leq \epsilon$ holds in expectation.
\label{cor:fixed}
\end{corollary}

\begin{remark} Let $\a_{\max}^2 = \max_i \| \A_i \|^2$ be the maximum squared row norm of $\A$. Given an initial error $\epsilon_0$ and a desired tolerance $\epsilon$ to the true solution, choosing the fixed step size 
$$ \alpha^* = \frac{p^3\epsilon \sigma_{\min}^2(\A)}{4(2-p)(1-p)\|\xopt\|^2 \sum_{i} \| \A_i \|^4 + 2p \epsilon \a_{\max}^2 \sigma_{\min}^2(\A) }, $$
after 
$$ k = 2 \log \left( \frac{2 \epsilon_0}{\epsilon} \right) \left(\frac{m\a_{\max}^2 }{p^2\sigma^2_{\min}(\A)} + \frac{2(2-p)(1-p) m \|\xopt \|^2\sum_i\| \A_i \|^4}{p^3 \sigma_{\min}^4(\A) \epsilon}\right)$$
iterations of Algorithm~\ref{alg:msgd}, $\E[ \| \xcur - \xopt \|^2 ]\leq \epsilon$ holds in expectation for consistent linear systems.
\end{remark}

\textbf{Recovering SGD. } When $p=1$, Algorithm~\ref{alg:msgd} behaves as classical SGD does on the full linear system $\A \x = \b$. Additionally, mSGD experiences similar convergence bounds as classical SGD for fixed step sizes \cite{needell2014stochastic}. In particular, when $p=1$ the updating function $g(\x)$ reduces to $g(\x) = \A_i^* (\A_i \x - \b_i) $.

\section{Experiments}
\label{sec:exp}
This section demonstrates the usefulness of Algorithm~\ref{alg:msgd} on synthetic and real world data. Although the full data set is available in every experiment, missing data is simulated by computing a binary mask that dictates which elements are available at every iteration. By doing so, the simplifying assumption is satisfied, and the ground truth is known, approximation error is computable, and we can investigate the performance of Algorithm~\ref{alg:msgd} with varying levels of missing data. In each experiment, the percentage of available data is varied and $\log$ $\ell_2$-error to the least squares solution, $\| \xcur - \xopt \|^2$ is averaged over 20 trials. For the fixed step size in simulated data,  $\alpha = 10^{-4}$ and for real world data  $\alpha = 10^{-5}$. For the updating step size regime, $\alpha_k =  \frac{c}{\sigma^2_{\min} (\A) k}$ with $c = 10^{-2}$. Using $\alpha_k = \frac{1}{\sigma^2_{\min} (\A) k}$ (as described in Theorem~\ref{thm:update}) creates an initial increase in error followed by a decrease in error. This behavior is attributed to the step sizes being too large initially. It seems that the factor $c$ can be optimized but we do not attempt to optimize such parameters here.

\begin{figure}[h!] 
    \centering
    \begin{subfigure}[b]{0.45\textwidth} 
        \includegraphics[width=\textwidth]{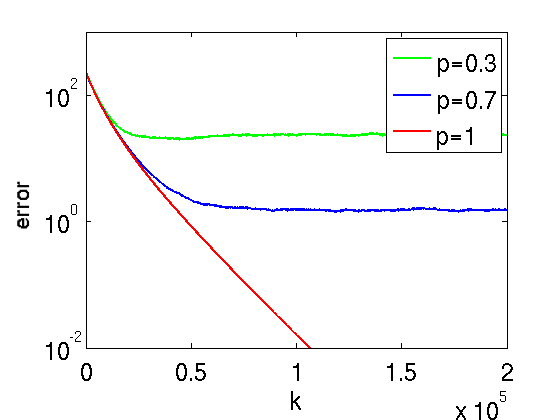}
        \caption{} \label{fig:gauss_consisa}
    \end{subfigure}
    ~
    \begin{subfigure}[b]{0.45\textwidth} 
        \includegraphics[width=\textwidth]{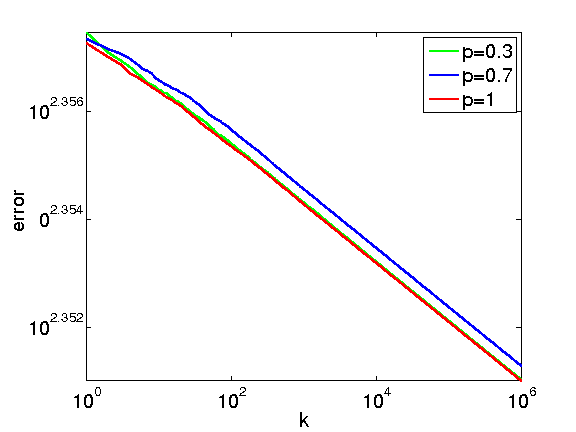}
        \caption{} \label{fig:gauss_consisb}
    \end{subfigure}
    \caption{This figure compares the performance of Algorithm~\ref{alg:msgd} on linear systems drawn from a standard Gaussian distribution. The percentage of data that is missing is varied. The x-axis is log(iteration) and the y-axis is the log($\ell_2$-error). Note that using a fixed step size (left), allows mSGD to converge much faster but to some convergence horizon. Using updating step sizes (right), continual progress is made at the cost of slower convergence. }  
    \label{fig:gauss_consis}
\end{figure}

\begin{figure}[h!]
    \centering
    \begin{subfigure}[b]{0.45\textwidth}  
        \includegraphics[width=\textwidth]{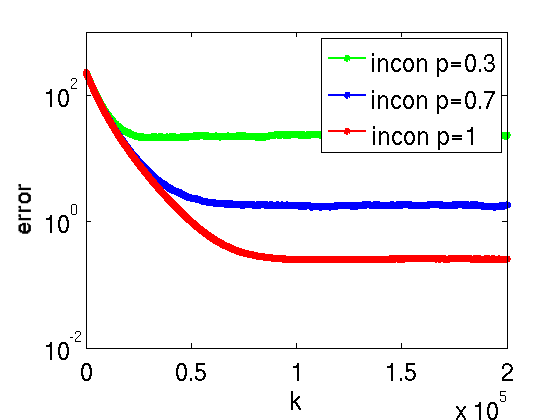}
        \caption{} \label{fig:gauss_incona}
    \end{subfigure}
    ~
    \begin{subfigure}[b]{0.45\textwidth}  
        \includegraphics[width=\textwidth]{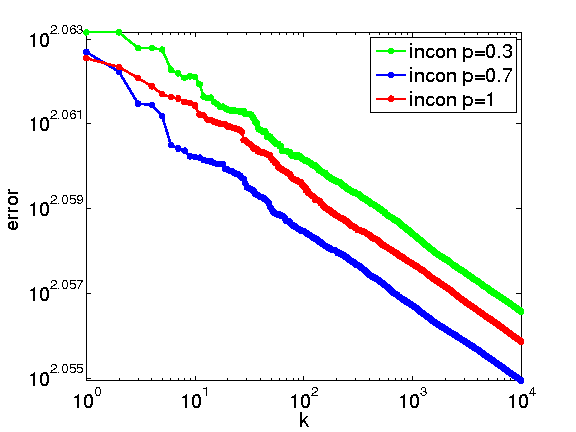}
        \caption{ } \label{fig:gauss_inconb}
    \end{subfigure} 
    \caption{The performance of mSGD on inconsistent linear systems. For a fixed step size (left), mSGD converges to a convergence horizon and using updating step sizes (right) allows mSGD to continually progress at a slower rate.}  \label{fig:gauss_incon}
\end{figure}

In the first experiment, we apply mSGD to synthetic data. The results can be seen in Figure~\ref{fig:gauss_consis} and Figure~\ref{fig:gauss_incon}. Here, elements of $\A \in \mathbb{R}^{m \times n}$ are drawn i.i.d. from a standard Gaussian distribution where $m = 1000 $ and $n = 200$. Figure~\ref{fig:gauss_consisa} and Figure~\ref{fig:gauss_incona} show the results of Algorithm~\ref{alg:msgd} using a fixed step size ($\alpha = 10^{-4}$) while Figure~\ref{fig:gauss_consisb} and Figure~\ref{fig:gauss_inconb} show results using updating step sizes. For inconsistent systems, we use $\b + \br$ as the right hand side vector where $\br$ is computed such that $\br \in \Null(\A^*)$ using Matlab's \texttt{null()} function.

\begin{figure}[h!] 
    \centering
    \begin{subfigure}[b]{0.45\textwidth}
        \includegraphics[width=\textwidth]{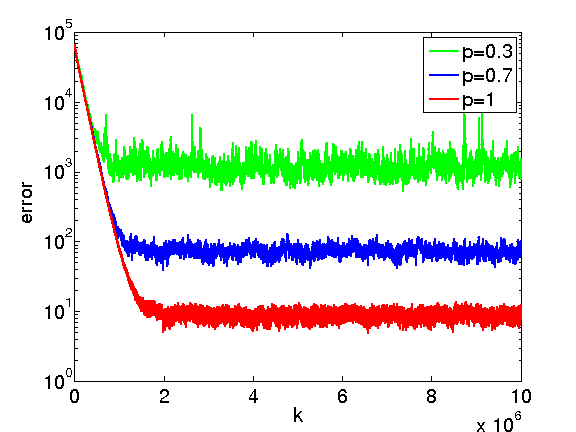}
        \caption{}
    \end{subfigure}
    ~
    \begin{subfigure}[b]{0.45\textwidth}
        \includegraphics[width=\textwidth]{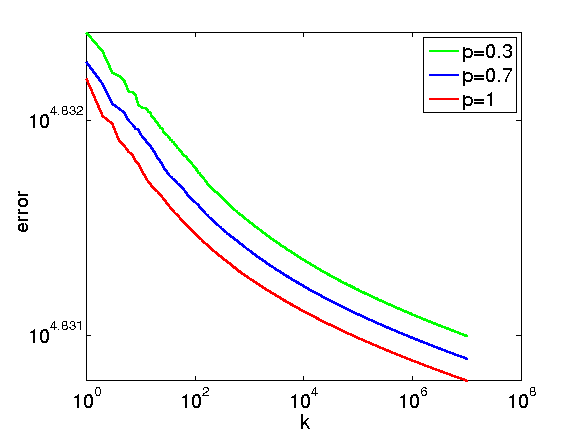}
        \caption{ }
    \end{subfigure}
    \caption{For bike data set, mSGD with a fixed step size (left) experiences a congerence horizon. Using updating step sizes (right), mSGD continues to progress toward the least squares solution.} \label{fig:bike}
\end{figure}

The first real world data set was obtained form the UCI Machine Learning Repository \cite{ucirepo} and contains data from a bike rental service. Rows of $\A$ contain hourly information from a bike share rental system and columns contain information such as weather, total number of rented bikes, time, and day of the week. In this experiment, $m = 17379$ and $n = 9$. Figure~\ref{fig:bike} displays the performance of mSGD on this data set for fixed and updating step sizes.

\begin{figure}[h!]
    \centering
        \includegraphics[width=.5\textwidth]{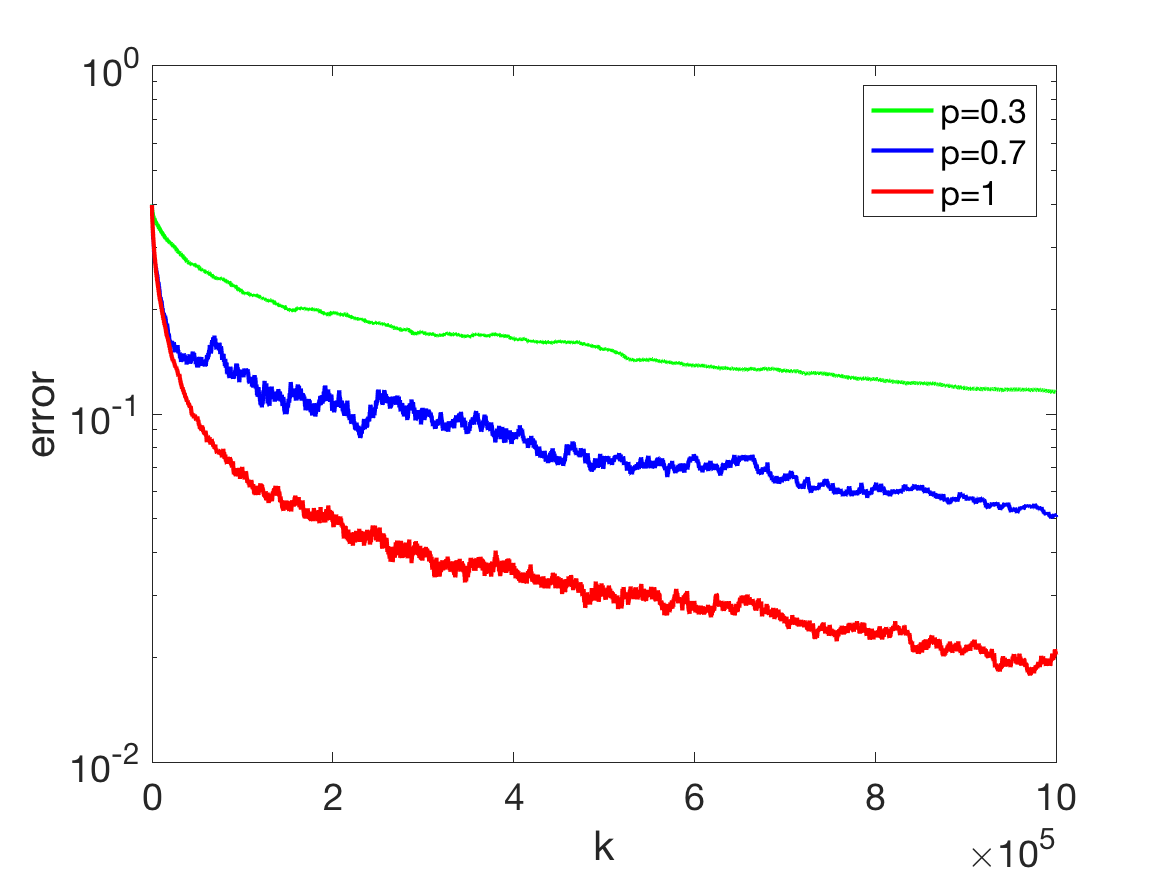}
        \caption{Using updating step sizes, Algorithm~\ref{alg:msgd} has decaying approximation error to the least squares solution of the completed linear system.}
        \label{fig:lyme}
\end{figure}

The performance of Algorithm~\ref{alg:msgd} on Lyme data from lymedisease.org is shown in Figure~\ref{fig:lyme}. This data set contains survey responses from patients who have been diagnosed with Lyme Disease. Examples of responses include number of emergency room visits, severity of symptoms, and effectiveness of medication.. For the right hand side vector, we use the number of health care providers a patient saw before being diagnosed with Lyme. For this experiment, $m=3686$ and $n = 81$. Solving such a system would potentially uncover what factors lead to late-stage diagnosis, a critical question in Lyme disease research. As seen in Figure~\ref{fig:gauss_inconb}, when using the updating step size $\alpha_k = \frac{c}{\sigma^2_{\min} (\A) k}$, the convergence rate suffers because the step size decays too quickly. Theoretically, we expect the error to continue to decay very slowly but, practically, it makes more sense to use another updating step size regime. In this experiment, we instead use $\alpha_k = \frac{c}{\sigma^2_{\min} (\A)} r^{\lfloor k/T* \rfloor}$ so that the initial step size is $\frac{c}{\sigma^2_{\min} (\A)}$ and after every $T^*$ iterations, the step size is multiplied by a factor of $r<1$. Empirical parameter tuning led us to use $c = 10^{-3}$ for $p = 0.7$ and $p=1$, $c=10^{-4}$ for $p=0.3$, $T^* = 10^5$, and $r = 0.8$. The results are shown in Figure~\ref{fig:lyme}.

\begin{figure}[h!] 
    \centering
        \includegraphics[width=.5\textwidth]{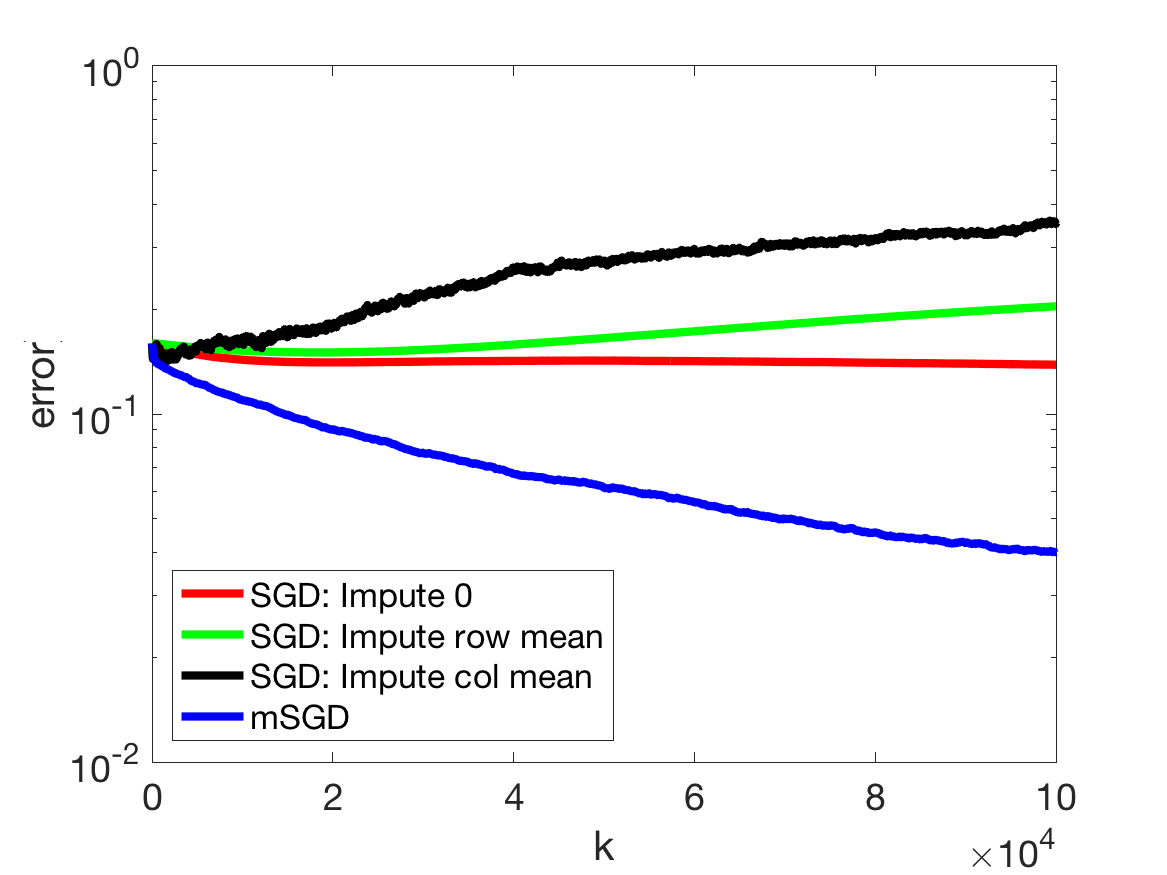}
  \caption{The proposed algorithm out performs using imputation methods with SGD.} \label{fig:impute}
\end{figure}

Figure~\ref{fig:impute} compares mSGD and classical SGD applied to three different imputation treatments for missing data. We use the Lyme Disease data set with updating step sizes as described in the previous experiment with $c=10^{-4}$ and $T^* = 10^5$. Setting $p=0.5$, $\A \in \mathbb{R}^{10^5 \times 81}$ where each row of $\A$ is a randomly selected row of the Lyme Disease data set with roughly half of the entries of the row (randomly) removed. Classical SGD is applied to $\Amiss$ in three ways: imputing 0 (if $\Amiss_{ij}$ is missing, $\Amiss_{ij} = 0$ ), imputing row means (if $\Amiss_{ij}$ is missing, $\Amiss_{ij}$ is the average over all non-missing elements in $\Amiss_i$), and imputing column means (if $\Amiss_{ij}$ is missing, $\Amiss_{ij}$ is the average over all non-missing elements in the $j^{th}$ column of $\Amiss$). Notice that mSGD outperforms the imputation methods presented here.

These experimental results support the theoretical findings presented in Section~\ref{sec:msgd}. Using a fixed step size, mSGD converges to some radius around the solution while using updating step size allows us to avoid the convergence horizon at the price of a slower convergence. For fixed step size, the amount of missing data affects the convergence horizon. In particular, as $p$ decreases the size of the convergence horizon increases.

\section{Conclusion}
\label{sec:conclu}
In this work, we present a stochastic iterative projection method that solves linear systems with missing data. We prove that mSGD finds the least squares solution to the linear system with full data even though a system has missing data. Additionally, this work shows theoretical bounds the performance of mSGD using fixed and updating step sizes. The experiments show that the proposed method is useful in real world settings when one wishes to solve a linear system with missing data without needing to impute missing values, which can be extremely costly.

\section{Appendix}

Consider the objective functions $F(\x) = \frac{1}{2m} \| \A \x - \b\|^2 = \frac{1}{m} \sum_{i=1}^m \frac{1}{2} (\A_i \x - \b_i)^2$ and $\tilde{F}(\x) = \frac{1}{2m} \| \Amiss \x - p\b \|^2$. Let $\tilde{f}_i(\x) = \frac{1}{2} (\amiss_i \x - p\b_i)^2$.  Let $\E_\delta [ \cdot ]$ denote the expected value function with respect to the Bernoulli random variables of the binary mask $\D$ and $\E_i [ \cdot]$ denote the expected value with respect to the choice of rows of $\Amiss$. In addition, let $\mu$ be the strong convexity parameter $F(\x)$ so that for any $\x, \y \in \mathcal{W}$, $(\x - \y)^*(\nabla F(\x) - \nabla F(\y) ) \geq \|\x - \y \|^2 \mu$.

First, we will show a few useful properties pertaining to the update function $g(\x)$. In particular, Lemma \ref{lem:egx} shows that in expectation $g(\x)$ allows us to make progress in the gradient direction of the objective $F(\x)$ (as opposed to the direction of $\nabla \tilde{F}(\x)$). Next, Lemma \ref{lem:lg} investigates the Lipschitz continuity of $g(\x)$ for a fixed row $i$ and binary mask $\D$ and its supremum Lipschitz constant of $g(\x)$ over all rows and binary masks. Lemma \ref{lem:G} shows that we can uniformly bound the expected norm of $g(\x)$ and provides said bound. Finally, we prove Theorem \ref{thm:fixed}.

\begin{lemma} The expected value of the update function $g(\x)$ defined in~\eqref{eq:update} is the gradient of the objection function $F(\x)$. In other words, for 
\begin{equation*}
g(\x) = \frac{1}{p^2} \left(\amiss_i^*(\amiss_i\x - p \b_i) \right) - \frac{(1-p)}{p^2}\diag(\amiss_i^*\amiss_i)\x,
\end{equation*}
we have that $\E [g(\x)] = \nabla F(\x)$.
\begin{proof}
To prove this lemma, we will first take the expected value of $\nabla \tilde{f}_i(\x)$. We then take the expected value of $g(\x)$, substitute $\E[\nabla \tilde{f}_i(\x)]$, and simplify to complete the proof. Let's first check that 
\begin{equation}
\E[\nabla \tilde{f}_i(\x)] = \E [\amiss_i^*(\amiss_i\x - p \b_i)] = p^2 \A^*\A\x + (p-p^2) \diag(\A^*\A)\x - p^2 \sum_i \ai \b_i.
\label{eq:egradfit}
\end{equation}
Taking a simple derivative, $\nabla \tilde{f}_i(\x) = \amiss_i^* (\amiss_i \x - p\b_i)$. The matrix $\D$ is a $m \times n$ binary mask with entries $\delta_{i,j} \stackrel{i.i.d.}{\sim} Bern(p)$. Let $\D_i = \diag(\delta_{i,1}, \delta_{i,2}, ... \delta_{i,n})$ be a $n \times n$ diagonal matrix so that $\amiss_i = \D_i \A_i^* $. Substituting $\amiss_i = \D_i \A_i^* $ and taking the expectation with respect to the $\delta_{i,j}$'s,  
\begin{align*}
\E_\delta [\nabla \tilde{f}_i(\x)] & = \E_\delta [\amiss_i^*\amiss_i ] \x - p\E_\delta[\amiss_i^*] \b_i \\
&=  \E_\delta [\amiss_i^*\amiss_i ] \x - p^2 \A_i^* \b_i \\
&\stackrel{(i)}{=} p^2 \A_i^*\A_i \x + (p-p^2) \diag(\A_i^*\A_i)\x - p^2 \A_i^* \b_i.
\end{align*}
Letting $[\A_i^* \A_i]_{jk}$ denote the $(j,k)^{th}$ element of $\A_i^*\A_i$, step ($i$) uses the fact that,
$$
 \E_\delta [\amiss_i^*\amiss_i ] = 
\begin{cases}
p [\A_i^* \A_i]_{jk}, &  j = k \\
p^2 [\A_i^* \A_i]_{jk}, &  j \neq k 
\end{cases}.
$$
Now, we take the expectation with respect to the rows of $\A$ to obtain:
\begin{align*}
\E [g(\x) ]
&\stackrel{(i)}{=}\frac{1}{p^2} \E [\amiss_i^*(\amiss_i\x - p \b_i) ] - \frac{1-p}{p^2} \E[\diag (\amiss_i^*\amiss_i)]\x \\
&=\frac{1}{p^2} \E [\nabla \tilde{f}_i(\x) ] - \frac{1-p}{p^2} \E[\diag (\amiss_i^*\amiss_i)]\x \\
&\stackrel{(ii)}{=}\frac{1}{ mp^2} \left( p^2 \A^*\A\x + (p-p^2)\diag(\A^*\A)\x - p^2 \sum_i \ai ^* \b_i \right) - \frac{p(1-p)}{ mp^2}\diag(\A^*\A)\x  \\
&=   \frac{1}{m}\A^*\A\x +\frac{(p-p^2)}{m p^2} \diag(\A^*\A)\x - \frac{1}{m} \A^*\b  - \frac{p-p^2}{m p^2}\diag(\A^*\A)\x  \\
&=  \frac{1}{m} \left(\A^*\A\x - \A^*\b \right) = \nabla F(\x).
\end{align*}
Step ($i$) follows from the definition of $g(\x)$ and linearity of the expected value. Step ($ii$) utilizes \eqref{eq:egradfit} for the first expected value and evaluates the expectation of $\E [\diag (\amiss_i^*\amiss_i)] = \E_i [\E_\delta [\diag (\amiss_i^*\amiss_i)] ]= p \E_i [\diag (\A_i^*\A_i)] = \frac{p}{m} \diag(\A^*\A)$. The remaining steps follow by simplification.
\end{proof}
\label{lem:egx}
\end{lemma}

\begin{lemma} The update function $g(\x)$ of Algorithm \ref{alg:msgd} is Lipschitz continuous with Lipschitz constant $L_{i,D}$. In other words, for all $\x,\y \in \mathcal{W}$, 
$$ \| g(\x) - g(\y) \| \leq L_{i,D} \| \x - \y \|. $$
In addition, we can bound the supremum Lipschitz constant, $L_g$  by
$$L_g = \sup_{i,D} L_{g,i,D} \leq \frac{a_{max}^2}{p^2},$$
where $a_{max}^2 = \max_i \| \A_i \|^2$. 
\begin{proof}
First we show that the Lipschitz constant $L_{i,D}$ of $g(\x)$
\begin{align*} 
\| g(\x) - g(\y) \| & =  \norm{ \left( \frac{1}{p^2} \amiss_i^*\amiss_i - \frac{(1-p)}{p^2} \diag(\amiss_i^*\amiss_i)\right)  \left(\x - \y\right)  }\\
& \leq \norm{ \frac{1}{p^2} \amiss_i\amiss_i^* - \frac{(1-p)}{p^2} \diag(\amiss_i^*\amiss_i) } \|\x - \y \| \\
& \leq \frac{1}{p^2} \norm{\amiss_i}^2  \|\x - \y \|.
\end{align*}
The last step follows from Weyl's Inequality which allows us to bound $\norm{ \amiss_i^*\amiss_i - (1-p) \diag(\amiss_i^*\amiss_i) } \leq \| \amiss_i \|^2$. 
Therefore we conclude that the Lipschitz constant of $g(\x)$ is $L_{i, D}  =  \frac{1}{p^2}\| \amiss_i\|^2$. 

To determine the supremum Lipschitz constant, we simply bound $L_{i,D}$ over all possible rows and all possible binary masks:
\begin{align*}
L_g & = \sup_{i,D} L_{i,D} = \sup_{i,D} \frac{1}{p^2}\norm{ \amiss_i}^2  \\
& \leq \frac{1}{p^2} \sup_i \norm{ \A_i}^2 \\
& \leq \frac{ a_{max}^2}{p^2} ,
\end{align*}
where $a_{max}^2$ is a largest row norm of $\A$.

\end{proof}
\label{lem:lg}
\end{lemma}

\begin{lemma} We can uniformly bound the expected value of the magnitude of the update function in the following way. We have that $E \| g(\x) \|^2 \leq G $, where
$$G = \frac{2B(2+p)(1-p)}{mp^3}  \sum_i \norm{ \A_i }^4  + \frac{2}{mp^2} \sum_i \norm{ \A_i }^2 | \b_i|^2,$$
where $B = \max_{\x \in \mathcal{W}} \| \x \|^2$. In addition, we have that 
\begin{itemize}
\item if $\A \xopt = \b$ (the linear system is consistent) then
$$ G_\star = \frac{2(1-p)(2 - p)}{mp^3}  \| \xopt \|^2 \sum_i \norm{ \A_i }^4.$$
\item if $\A \xopt = \b + \br$ (the linear system is inconsistent) then
$$ G_\star = \frac{2}{mp^2} \sum_i \norm{ \A_i }^2  \br_i^2  + \frac{2(1-p)(2 - p)}{mp^3}  \| \xopt \|^2 \sum_i \norm{ \A_i }^4.$$
\end{itemize}
$G$ and $G_\star$ are also defined in \eqref{eq:G} and \eqref{eq:Gstar} respectively.
\begin{proof}

We begin this proof by showing the upper bound of $\E \left[ \|g(\x) \|^2 \right]$ for a general $\x$. From here, we obtain $G_\star$ by substituting $\x$ with $\xopt$ and making the appropriate assumptions on the consistency of the linear system. To get the uniform upper bound over all $\x$, we isolate $\| \x \|^2$ and bound the norm by $B = \max_{\x \in \mathcal{W}} \|  \x \|^2$. We have,

\begin{align*}
\E [\| g(\x) \|^2] & =  \E \left[ \norm{ \frac{1}{p^2} \left(\amiss_i^*(\amiss_i\x - p \b_i) \right) - \frac{(1-p)}{p^2}\diag(\amiss_i^*\amiss_i)\x }^2 \right] \\ 
&\stackrel{(i)}{\leq} \frac{2}{p^4} \E \left[\norm{ \amiss_i^*(\amiss_i\x - p \b_i) }^2 \right]+ \frac{2(1-p)^2}{p^4} \E \left[\| \diag(\amiss_i^*\amiss_i)\x \|^2\right] \\ 
&\stackrel{(ii)}{=} \frac{2}{p^4} \E \left[\| \amiss_i \|^2 (\amiss_i\x - p \b_i)^2\right] + \frac{2(1-p)^2}{p^4} \E \left[\| \diag(\amiss_i^*\amiss_i)\x \|^2\right] \\ 
& \stackrel{(iii)}{\leq} \frac{2}{p^4} \E \left[\norm{ \A_i }^2 (\amiss_i\x - p \b_i)^2\right] + \frac{2(1-p)^2}{p^4} \E \left[ \| \diag(\amiss_i^*\amiss_i)\x \|^2\right].
\end{align*}

Step ($i$) follows by Jensen's inequality, step ($ii$) is simplification and uses the fact that $(\amiss_i\x - p \b_i)$ is scalar. Lastly, step ($iii$) bounds the magnitude of a row of $\A$ with missing data by the magnitude of a row of $\A$ without missing data (i.e. $\| \tilde{\A}_i \| = \| \D_i \A_i \| \leq \| \A_i \|$ for all $\D_i$). From here, we use the fact that $\E = \E_i \E_\delta$ to obtain the following: 

\begin{equation}
\E[ \| g(\x) \|^2 ]\leq \frac{2}{p^4} \E_i \big [\norm{ \A_i }^2 \underbrace{\E_\delta [(\amiss_i\x - p \b_i)^2]}_\text{(A)} \big ]+ \frac{2(1-p)^2}{p^4} \E_i \bigg [\underbrace{\E_\delta \left[ \| \diag(\amiss_i^*\amiss_i)\x \|^2 \right]}_\text{(B)} \bigg ].
\label{eq:mid}
\end{equation}
Now, we will focus on the computation of $\E_\delta$. First, we compute (A). We have that,

\begin{align*}
\E_\delta& \left[ (\amiss_i\x - p \b_i)^2 \right]
 =  \E_\delta \left[ (\amiss_i\x )^2 \right] - 2p \E_\delta \left[\amiss_i\right]\x \b_i  + p^2 \b_i^2 \\
& =  \E_\delta \left[ \left(\sum_{j=1}^n \Amiss_{ij} \x_j \right)^2 \right] - 2p^2 \A_i\x \b_i + p^2 \b_i^2 \\
& = \E_\delta \left[ \sum_{j=1}^n \Amiss_{ij}^2 \x_j ^2 + 2 \sum_{j=1}^n \sum_{k=1}^{j-1}  \Amiss_{ij}\Amiss_{ik} \x_j \x_k \right] - 2p^2  \A_i\x \b_i + p^2 \b_i^2 \\
& = \left( p \sum_{j=1}^n \A_{ij}^2 \x_j ^2 + 2p^2 \sum_{j=1}^n \sum_{k=1}^{j-1}  \A_{ij}\A_{ik} \x_j \x_k \right) - 2p^2 \A_i\x\b_i  + p^2 \b_i^2  \\
& \stackrel{(i)}{=} \left( p^2 \sum_{j=1}^n \A_{ij}^2 \x_j ^2 + (p-p^2) \sum_{j=1}^n \A_{ij}^2 \x_j ^2 + 2p^2 \sum_{j=1}^n \sum_{k=1}^{j-1}  \A_{i,j}\A_{i,k} \x_j \x_k \right) - 2p^2 \A_i\x \b_i  + p^2 \b_i^2  \\
& =  p^2 \left( \sum_{j=1}^n \A_{ij}^2 \x_j ^2 + 2\sum_{j=1}^n \sum_{k=1}^{j-1}  \A_{ij}\A_{ik} \x_j \x_k - 2\A_i\x\b_i + \b_i^2 \right)  + (p-p^2) \left( \sum_{j=1}^n \A_{ij}^2 \x_j ^2\right) \\
& =  p^2 \left( \left(\sum_{j=1}^n \A_{ij} \x_j\right)^2 - 2\A_i\x\b_i  + \b_i^2 \right)  + (p-p^2) \x^* \diag(\A_i^*\A_i) \x\\
& =  p^2 \left( \A_i \x - \b_i \right)^2  + p(1-p) \x^* \diag(\A_i^*\A_i) \x .
\end{align*}

In step ($i$), we add and subtract the term $p^2 \sum_{j=1}^n \A_{i,j}^2\x_j^2$ so that we can combine terms. Other equalities follow by simplification and computation of expected value. Note that $\E_\delta [ \tilde{\A}_{i,j} ] = p\A_{i,j}$ and $\E_\delta [ \tilde{\A}_{i,j}\tilde{\A}_{i,k}] = p^2\A_{i,j}\A_{i,k}$ if $j \neq k$. 

For term (B), we simply compute that 
\begin{align*}
\E_\delta \left[ \| \diag(\amiss_i^*\amiss_i)\x \|^2 \right] &= \E_\delta \left[ \sum_{j=1}^n \tilde{\A}_{i,j}^2 \x_j^2 \right] \\
& = p \sum_{j=1}^n (\A_{i,j}^2 \x_j^2) \\
&= p \norm{ \diag(\A_i^*\A_i)\x }^2.
\end{align*}
 Now that we have (A) and (B), we can compute a general upper bound for $\E [\|g(\x) \|^2]$. Starting with substituting (A) and (B) into \eqref{eq:mid},

\begin{align*}
\E &\left[ \| g(\x) \|^2 \right] \stackrel{(i)}{\leq} \frac{2}{p^2} \E_i \left[ \norm{ \A_i }^2  \left(\A_i \x - \b_i \right)^2 \right]   \\
& \hspace{2cm}+ \frac{2p(1-p)}{p^4} \E_i \left[ \norm{ \A_i }^2 \x^* \diag(\A_i^*\A_i) \x \right] + \frac{2p(1-p)^2}{p^4} \E_i \left[ \norm{ \diag(\A_i^*\A_i) \x }^2 \right] \\ 
&\hspace{1cm}\stackrel{(ii)}{\leq} \frac{2}{p^2} \E_i \left[ \norm{ \A_i }^2  \left( \A_i \x - \b_i \right)^2 \right] + \left( \frac{2p(1-p)}{p^4} +  \frac{2p(1-p)^2}{p^4} \right) \E_i \left[ \norm{ \A_i }^2 \x^* \diag(\A_i^*\A_i) \x  \right] \\
&\hspace{1cm}\leq \frac{2}{p^2} \E_i \left[ \norm{ \A_i }^2  \left( \A_i \x - \b_i \right)^2  \right] + \frac{2p(1-p)(2 - p)}{p^4} \E_i  \left[ \norm{ \A_i}^2\x^*\diag(\A_i^*\A_i) \x \right] \\
&\hspace{1cm}\stackrel{(iii)}{=} \frac{2}{mp^2} \sum_i \norm{ \A_i }^2  \left( \A_i \x - \b_i \right)^2  + \frac{2p(1-p)(2 - p)}{mp^4}  \sum_i \norm{ \A_i}^2\x^* \diag(\A_i^*\A_i) \x \\
&\hspace{1cm}\stackrel{(iv)}{\leq} \frac{2}{mp^2} \sum_i \norm{ \A_i }^2  \left( \A_i \x - \b_i \right)^2  + \frac{2p(1-p)(2 - p)}{mp^4}  \| \x \|^2 \sum_i \norm{ \A_i }^4. 
\end{align*}
Step ($i$) substitutes (A) and (B) in \eqref{eq:mid}. Step ($ii$) uses the fact that $\norm{ \diag(\ai^*\ai) \x }^2 \leq \norm{ \diag(\ai ) \diag(\ai ) \x}^2 \leq \norm{\ai}^4 \| \x\|^2$.

From here, we substitute $\x$ with $\xopt$ to compute $G_\star$. If $\A \xopt = \b$ (the linear system is consistent) then the terms $\left( \A_i \x - \b_i \right)^2 = 0$ and we find that 
$$G_\star = \frac{2(1-p)(2 - p)}{mp^3}  \| \xopt \|^2 \sum_i \norm{ \A_i }^4.$$

 Otherwise, if $\A \xopt = \b+\br$ for some residual vector $\br$, we have that 
 $$G_\star = \frac{2}{mp^2} \sum_i \norm{ \A_i }^2  \br_i^2  + \frac{2(1-p)(2 - p)}{mp^3}  \| \xopt \|^2 \sum_i \norm{ \A_i }^4,$$ 
 where $\br_i$ is the $i^{th}$ element of the vector $\br$. To finish the proof of Lemma \ref{lem:G}, we simplify starting from step ($iv$).

\begin{align*}
\E \left[ \| g(\x) \|^2 \right]  &\leq \frac{2}{mp^2} \sum_i \norm{ \A_i }^2  \left( \A_i \x - \b_i \right)^2  + \frac{2p(1-p)(2 - p)}{mp^4}  \| \x \|^2 \sum_i \norm{ \A_i }^4 \\
& \stackrel{(i)}{\leq} \frac{2}{mp^2} \sum_i \norm{ \A_i }^2  \left( | \A_i \x |^2 + | \b_i|^2 \right)  + \frac{2p(1-p)(2 - p)}{mp^4}  \| \x \|^2 \sum_i \norm{ \A_i }^4 \\
& \stackrel{(ii)}{\leq} \frac{2}{mp^2} \sum_i \norm{ \A_i }^4   \| \x \|^2 + \frac{2}{mp^2} \sum_i \norm{ \A_i }^2 | \b_i |^2   + \frac{2p(1-p)(2 - p)}{mp^4}  \| \x \|^2 \sum_i \norm{ \A_i }^4 \\
& \stackrel{(iii)}{\leq} \frac{2}{mp^2} \sum_i \norm{ \A_i }^4  B + \frac{2}{mp^2} \sum_i \norm{ \A_i }^2 | \b_i |^2   + \frac{2p(1-p)(2 - p)}{mp^4} B \sum_i \norm{ \A_i }^4 \\
& = \left( \frac{2B}{mp^2} + \frac{2p(1-p)(2 - p)B}{mp^4}  \right) \sum_i \norm{ \A_i }^4  + \frac{2}{mp^2} \sum_i \norm{ \A_i }^2 | \b_i|^2 
\\
& = \frac{2B}{mp^2}\left( 1 + \frac{p(1-p)(2 - p)}{p^2}  \right) \sum_i \norm{ \A_i }^4  + \frac{2}{mp^2} \sum_i \norm{ \A_i }^2 | \b_i |^2 
\\
& = \frac{2B}{mp^2}\left( 1 + \frac{(1-p)(2 - p)}{p}  \right)  \sum_i \norm{ \A_i }^4  + \frac{2}{mp^2} \sum_i \norm{ \A_i }^2 | \b_i|^2 
\end{align*}
In step ($i$) we use Jensen's inequality. Note that $\A_i \x$ and $\b_i$ are both scalar values. In step ($ii$), we distribution the summation in the first term and use the fact that $|\A_i \x |^2 \leq \| \A_i \|^2 \| \x \|^2$ by the Cauchy-Schwarz inequality. Step ($iii$) uses the definition of $B = \max_{\x \in \mathcal{W}} \| \x \|^2$. The remaining lines are simplification.
\end{proof}
\label{lem:G}
\end{lemma}

Before we begin the proof of Theorem \ref{thm:fixed}, we remind the reader that $F(\x)$ is strongly convex with strong convexity parameter $\mu$. In other words, for all $\x, \y \in \mathcal{W}$ we have that
\begin{equation}
(\x - \y)^*(\nabla F(\x) - \nabla F(\y)) \geq \mu \|\x - \y \|^2.
\label{eq:strongcon}
\end{equation}
In addition, we define a new function $G(\x) = \frac{1}{2p^2} \left( (\tilde{\A}_i \x - p\b_i )^2 - \frac{(1-p)}{2p^2} \|\diag(\tilde{\A}_i) \x \|^2\right) $ so that $g(\x) = \nabla G(\x)$. The update function $g(\x)$ follows the Co-coercivity Lemma as stated in Lemma \ref{lem:cocoer}.

\begin{lemma}(\cite{needell2014stochastic} Lemma A.1) For $G(\x)$ a smooth function such that $\nabla G(\x) = g(\x)$,
\begin{equation*}
\|g(\x) - g(\y) \|^2 \leq L_{i,D} (\x - \y)^*(g(\x) - g(\y)),
\end{equation*}
where $g(\x)$ has Lipschitz constant $L_{i,D}$.
\label{lem:cocoer}
\end{lemma}

\subsection{Proof of Theorem \ref{thm:fixed}}

\begin{proof}
First, we bound expected error conditional on the previous $k-1$ iterations. Let $\E_{k-1} [\cdot]$ denote the expected value conditional of the previous $k-1$ iterations and note that by the Law of Iterated Expectation, we have that the full expected value over all iterations is $\E[\cdot] = \E[\E_{k-1}[\cdot]] $. Thus,
	\begin{align*}
	\E_{k-1} \| \xcur - &\xopt \|^2 = \E_{k-1} \left[ \| \xprev - \alpha g(\xprev) - \xopt \|^2 \right] \\
	& \stackrel{(i)}{=} \| \xprev  - \xopt \|^2 - 2 \alpha (\xprev - \xopt)^* \E_{k-1} \left[ g(\xprev) \right] +  \alpha^2 \E_{k-1} \left[ \| g(\xprev) \|^2 \right]  \\
	& \stackrel{(ii)}{=} \| \xprev  - \xopt \|^2 - 2 \alpha (\xprev - \xopt)^* (\nabla F(\xprev) - \nabla F(\xopt)) +  \alpha^2 \E_{k-1} \left[ \| g(\xprev) \|^2 \right] \\
	& \stackrel{(iii)}{\leq} \| \xprev  - \xopt \|^2 - 2 \alpha (\xprev - \xopt)^* (\nabla F(\xprev) - \nabla F(\xopt)) \\
	    &  \hspace{5em} + 2 \alpha^2 \E_{k-1} \left \| g(\xprev) - g(\xopt) \|^2 \right] +  \alpha^2\E_{k-1} \left[ \| g(\xopt)\|^2 \right] \\
	& \stackrel{(iv)}{\leq} \| \xprev  - \xopt \|^2 - 2 \alpha (\xprev - \xopt)^* (\nabla F(\xprev) - \nabla F(\xopt)) \\
	     &  \hspace{5em} + 2\alpha^2 L_{i,g} (\xprev - \xopt)^*(\E_{k-1} [ g(\xprev)] -  \E_{k-1}  [g(\xopt)])  +  \alpha^2 G_\star \\
	& \stackrel{(v)}{\leq}  \| \xprev  - \xopt \|^2 - 2 \alpha (\xprev - \xopt)^* (\nabla F(\xprev) - \nabla F(\xopt))  \\
	      & \hspace{5em}+ 2 \alpha^2 L_g (\xprev - \xopt)^*(\nabla F(\xprev) -  \nabla F(\xopt)) + \alpha^2 G_\star \\
	& \leq  \| \xprev  - \xopt \|^2 - 2 \alpha (1- \alpha L_g) (\xprev - \xopt)^* (\nabla F(\xprev) - \nabla F(\xopt))  + \alpha^2 G_\star \\
		& \stackrel{(vi)}{\leq} \| \xprev  - \xopt \|^2 - 2 \alpha (1- \alpha L_g) \mu \|\xprev - \xopt \|^2  + \alpha^2 G_\star \\ 
		& = \left(1 - 2 \alpha \mu (1- \alpha L_g) \right) \| \xprev - \xopt \|^2 + \alpha^2 G_\star \\
		& = r \| \xprev - \xopt \|^2 + \alpha^2 G_\star .
\end{align*}
Step ($i$) follows from the definition of the $\ell_2$ norm. Step ($ii$) takes the expectation of $g(\xprev)$ using Lemma~\ref{lem:egx} and uses the fact that $\nabla F(\xopt) = 0$ to subtract $2\alpha(\xprev - \xopt)^*\nabla F(\xopt)$. In step ($iii$) we add and subtract the term $\|g(\xopt) \|^2$ then apply Jensen's inequality. Step ($iv$) is an application of the Lemma \ref{lem:cocoer}. Step ($v$) bounds $L_{i,D}$ by $L_g = \sup_{i,D} L_{i,D}$ and uses Lemma \ref{lem:egx} to compute the expectation of $\E_{k-1} [ g(\x)]$. We use the strong convexity of $F(\x)$ in step ($vi$). The remaining lines are simplification. Now, by the Law of Iterated Expectation we recursively apply this bound to obtain the desired result,

\begin{align*}
\E \| \xcur - \xopt \|^2 & \leq  r \E_{k-2} \| \xprev  - \xopt \|^2 + \alpha^2 G_\star	\\
& \leq r^k \| \x_0 - \xopt \|^2 + \alpha^2 G_\star \sum_{j = 0}^{k-1} r^j  \\
& \leq r^k \| \x_0 - \xopt \|^2 + \frac{\alpha^2 G_\star }{1- r}.
\end{align*}

\textbf{A note on Inconsistent Linear Systems. } Theorem \ref{thm:fixed} also applies to inconsistent systems. Let $\A \xopt = \b + \br$ where $\br \in \Null(\A^*)$. In the proof of Theorem \ref{thm:fixed}, we use the fact that $\nabla F(\xopt) = 0$ in step ($ii$). This is still true in the inconsistent setting as $\nabla F(\xopt) = \A^*(\A \xopt - \b) = \A^* \br = 0$. All other computations go through without issue. 

\end{proof}

\section*{Acknowledgments}
Needell was partially supported by NSF CAREER grant $\#1348721$, NSF BIGDATA $\#1740325$,  and the Alfred P. Sloan Fellowship. Ma was supported in part by NSF CAREER grant $\#1348721$, the CSRC Intellisis Fellowship, and the Edison International Scholarship.

\bibliographystyle{abbrv}
\bibliography{bib_am}
\end{document}